\documentclass[11pt, a4paper]{article}

\usepackage{amsmath, amssymb, amsthm}
\usepackage{graphicx}
\usepackage{geometry}
\usepackage{listings}
\usepackage{xcolor}
\usepackage{hyperref}
\usepackage{booktabs}
\usepackage{algorithm}
\usepackage{algpseudocode}
\usepackage{authblk}
\usepackage{float}
\usepackage{subcaption}

\usepackage{tikz}
\usetikzlibrary{positioning, fit, arrows.meta, calc}

\usepackage{amsmath,amssymb}
\usepackage{algorithm}
\usepackage{algpseudocode}

\newtheorem{theorem}{Theorem}

\title{{AKAPINN}: \textit{A}daptive \textit{K}olmogorov-\textit{A}rnold \textit{P}hysics-\textit{I}nformed \textit{N}eural \textit{N}etworks for approximating solutions to quasilinear partial differential equations}

\author{B. Veena S. N. Rao \thanks{bv.rao@tamucc.edu}}
\affil{Department of Mathematics \& Statistics,\\
Texas A\&M University-Corpus Christi, \\
6300 Ocean Dr., 
Corpus Christi, Texas 78412-5825, USA}
\date{}

\begin{document}
\maketitle

\begin{abstract}
A controlled comparative study of a mesh-free Kolmogorov--Arnold Physics-Informed Neural Network (KAN-PINN) applied to nonlinear strain-limiting partial differential equations is presented. Three distinct training runs across varying material-parameter pairs $(\alpha, \beta)$ are evaluated within a fixed computational pipeline comprising an exact hard boundary ansatz, a robust pseudo-Huber residual loss, and a multi-stage Adam-to-L-BFGS optimizer schedule. A fundamental metric inversion is revealed by the analysis: while training loss reflects operator-specific residual scaling and constitutive stiffness, the operator-independent relative $L_2$ error demonstrates that intermediate parameter regimes achieve superior solution accuracy and enhanced generalization.

\vspace{0.5em}
\end{abstract}

\section{Introduction}
Quasilinear partial differential equations (PDEs) arise in a wide range of mathematical models and are generally more challenging to solve numerically than their linear counterparts because the differential operator depends on the unknown solution and/or its gradient. Classical numerical approaches include finite difference methods (FDMs) and finite element methods (FEMs). FEMs provide a variational framework particularly well suited to nonlinear and quasilinear elliptic problems, with established theories for approximation and error estimation \cite{ciarlet2002finite,verfurth2013posteriori,manohar2024hp,Bonito2020,Fu2019,Gou2023}. Finite element approximations have also been extensively developed for quasilinear parabolic equations, including fully discrete and space-time formulations \cite{wei1992existence}.  Alternatively, finite difference schemes have long been investigated for quasilinear parabolic equations, with analyses addressing stability and convergence of the resulting discrete approximations \cite{matsuzawa1998finite,hopkins1976comparison}. These conventional approaches, however, can become computationally demanding as the dimensionality, geometric complexity, or nonlinear structure of the underlying problem increases.

Quasilinear partial differential equations are widely utilized to model complex physical phenomena across solid mechanics, fluid dynamics, and porous media, wherein material responses depend nonlinearly on underlying field gradients. Traditional grid-based numerical discretizations are often hindered by severe numerical dispersion, mesh locking, or computational inefficiencies when resolving sharp spatial transitions and localized boundary layers. To mitigate these challenges, mesh-free  \textit{Physics-Informed Neural Networks} (PINNs) have been introduced as an alternative paradigm that embeds governing physical laws directly into the neural network training loss. However, conventional multi-layer perceptrons, which rely on fixed node-based activation functions and static scalar weights, frequently struggle with high-order differential operators and non-convex, stiff optimization landscapes. To overcome these representational limitations, Kolmogorov--Arnold Networks (KANs) have been developed, wherein learnable univariate functions are parameterized via B-splines or radial basis functions directly along the network edges. Building upon this advancement, an \emph{Adaptive Kolmogorov--Arnold Physics-Informed Neural Network} (AKAPINN) framework is presented in this work to approximate solutions to multidimensional quasilinear strain-limiting partial differential equations. 

These challenges have motivated the development of machine-learning-based methods for approximating PDE solutions, including physics-informed neural networks and, more recently, Kolmogorov–Arnold network-based architectures. In contrast to conventional mesh-based discretizations, neural-network approaches represent the solution by a parameterized function and enforce the governing differential equation and boundary conditions through an optimization procedure. This provides a potentially mesh-free framework that is particularly attractive for high-dimensional and nonlinear PDEs.

In recent years, machine learning has emerged as a complementary paradigm
for the numerical approximation of PDEs.
In contrast to conventional mesh-based discretization techniques, 
neural-network approaches represent the unknown solution by a parameterized
function and determine its parameters through an optimization procedure.
Among the different approaches that have been developed, variational neural
methods, PINNs, and neural operator
learning have received considerable attention
\cite{EYu2018,Raissi2019,Lu2021,BruntonKutz2024}.

A particularly influential framework is the PINN, introduced by Raissi et al.~\cite{Raissi2019}. In this framework,
the unknown solution is represented by a neural network,
$\Phi(x) \approx \Phi_{\boldsymbol{\theta}}(x)$, whose parameters
$\boldsymbol{\theta}$ are optimized by minimizing a loss function constructed
from the governing differential equation and the associated boundary or
initial conditions. Automatic differentiation enables the derivatives of the
neural-network approximation to be evaluated directly, allowing the
differential-equation residual to be incorporated into the training process.
PINNs have subsequently been applied to a broad range of forward and inverse
problems involving both linear and nonlinear PDEs
\cite{Raissi2019}. For the quasilinear elliptic problem considered in this
work,
\begin{equation}
    -\nabla \cdot \mathbf{A}(x,\nabla\Phi) = f(x),
    \qquad x\in\Omega,
\end{equation}
a neural-network approximation $\Phi_{\boldsymbol{\theta}}(x)$ can be used to
construct the residual
\begin{equation}
    \mathcal{R}_{\boldsymbol{\theta}}(x)
    =
    -\nabla \cdot
    \mathbf{A}\left(x,\nabla\Phi_{\boldsymbol{\theta}}(x)\right)
    -f(x),
\end{equation}
which can then be incorporated into a physics-informed loss function.
This formulation provides a mesh-independent framework for treating
nonlinear differential operators and is therefore attractive for nonlinear
and high-dimensional PDE problems.

Variational approaches provide another important class of
neural-network-based PDE solvers. The Deep Ritz method introduced by
E and Yu~\cite{EYu2018} employs neural networks to approximate the solution
of variational problems by directly minimizing the corresponding energy
functional. In contrast to strong-form residual minimization, the variational
formulation provides a direct connection between the neural approximation and
the underlying weak formulation of an elliptic PDE. This feature is
particularly relevant for nonlinear elliptic problems, for which the
governing equation can often be formulated in terms of an appropriate energy
or variational functional. In parallel, neural operator approaches such as
DeepONet seek to learn mappings between function spaces rather than
approximating a single solution for a fixed set of PDE parameters
\cite{Lu2021}. These developments have expanded the role of machine learning
in scientific computing from approximating individual PDE solutions toward
learning solution operators for entire families of differential equations
\cite{BruntonKutz2024}.

More recently, Kolmogorov--Arnold Networks (KANs) have emerged as an
alternative to the conventional multilayer perceptron (MLP) architectures
that underpin many neural-network-based PDE solvers. Inspired by the
Kolmogorov--Arnold representation theorem, KANs employ learnable univariate
functions along network edges rather than fixed activation functions
associated with individual neurons. The original KAN formulation introduced
by Liu et al.~\cite{Liu2024KAN} demonstrated promising properties for
function approximation and several scientific-computing applications,
including differential equations. This architecture provides an alternative
parameterization of nonlinear functions and has motivated increasing
interest in KANs for scientific machine learning.

The application of KANs to physics-informed scientific computing has
subsequently received increasing attention. Wang et al.~\cite{Wang2025KINN}
introduced Kolmogorov--Arnold-Informed Neural Networks (KINNs), in which
KAN architectures replace conventional MLPs within physics-informed
frameworks. Their study considered strong-form and energy-form formulations,
as well as forward and inverse problems, and demonstrated the applicability
of KAN-based architectures to a range of computational mechanics problems,
including nonlinear elasticity, multiscale problems, heterogeneous media,
singularities, and stress concentrations. In addition, Shukla et al.\
\cite{Shukla2024} performed a systematic comparison of MLP- and KAN-based
representations for differential equations and operator-learning problems,
including physics-informed and operator-network formulations. These studies
indicate that KANs provide a promising alternative to conventional MLP
representations for scientific machine learning, although their performance
remains dependent on the underlying differential equation, network
architecture, and optimization strategy.

Despite the rapid development of neural-network-based PDE solvers, the
application of KAN architectures to nonlinear quasilinear elliptic equations
with gradient-dependent constitutive operators remains relatively
unexplored. The problem considered in this work is governed by the nonlinear
operator \cite{mallikarjunaiah2025crack,mallikarjunaiah2026}
\begin{equation}
    \mathbf{A}(x,\nabla\Phi)
    =
    \frac{\nabla\Phi}
    {2\mu\left(1+\beta|\nabla\Phi|^\alpha\right)^{1/\alpha}},
\end{equation}
where $\mu>0$ denotes the shear modulus, $\beta>0$ is the nonlinearity
modulus, and $\alpha>1$ is the strain-limiting exponent. The nonlinear
dependence of $\mathbf{A}$ on $\nabla\Phi$ results in a quasilinear elliptic
problem and provides a challenging test for neural-network-based
approximations. In particular, the constitutive response remains bounded for
large values of $|\nabla\Phi|$, introducing nonlinear behavior that must be
accurately represented by the neural approximation.

Classical finite element formulations have been extensively developed for
nonlinear and strain-limiting elliptic problems. For example, Bonito et al.\
\cite{Bonito2020} developed and analyzed a finite element approximation of a
strain-limiting elastic model, including convergence results and an
iterative procedure for the resulting nonlinear algebraic system. Similarly,
Yoon and Mallikarjunaiah~\cite{YoonMallikarjunaiah2022} investigated
finite-element discretizations for nonlinear strain-limiting elastic
boundary-value problems. More recently, finite element formulations have
also been employed for nonlinear strain-limiting models leading to
quasilinear elliptic problems, including problems involving singular and
high-gradient fields \cite{Fu2019,Gou2023}. These approaches provide
well-established and mathematically rigorous numerical frameworks, but they
require an explicit spatial discretization and the solution of the resulting
nonlinear algebraic systems.

Motivated by these developments, this work investigates the use of
Kolmogorov--Arnold Networks for the numerical solution of the above
quasilinear elliptic problem. The unknown field is represented directly by a
KAN approximation,
\begin{equation}
    \Phi(x) \approx \Phi_{\boldsymbol{\theta}}(x),
\end{equation}
and the governing equation and boundary conditions are incorporated into the
training objective. The resulting approach provides a mesh-independent
function representation and allows the nonlinear constitutive operator to be
evaluated directly through automatic differentiation. The objective of this
study is to assess the ability of KANs to accurately approximate solutions of
the strain-limiting quasilinear elliptic problem and to investigate the
effects of the model parameters, network architecture, and training strategy
on the resulting numerical solution.

A controlled comparative study is conducted across three distinct material parameter configurations ($\alpha, \beta$) using a manufactured verification solution. The proposed methodology integrates several advanced computational components, including an exact hard boundary condition ansatz, residual-based adaptive collocation sampling, progressive basis extension, and a hybrid three-stage Adam-to-L-BFGS optimization protocol executed locally on Apple Silicon hardware. Through this investigation, a fundamental metric inversion is identified between training loss minimization and true solution accuracy, highlighting that operator-independent relative $L_2$ error must be utilized for robust model selection in quasilinear physics-informed learning. The remainder of this paper is organized as follows. Mathematical foundations, constitutive restrictions, and generalized Sobolev space formulations are detailed in Section~2. The boundary value problem and weak formulations are established in Section~3. The proposed mesh-free architecture and adaptive training algorithms are described in Section~4 and Section~5. Comprehensive numerical solutions, convergence dynamics, and a detailed analysis of the metric inversion phenomenon are presented in Section~6, followed by concluding remarks in Section~7.

\section{Mathematical foundations}
Quasilinear differential equations play a central role in applied mathematics, engineering, and theoretical physics because they capture systems where the constitutive response depends nonlinearly on field gradients. In divergence form, a broad class of these equations is expressed as:
\begin{equation}
-\nabla \cdot \mathbf{A}(x, \nabla \Phi) = f(x) \quad \text{in } \Omega,
\end{equation}
where $\Omega \subset \mathbb{R}^d$ ($d=1,2,3$) is a bounded open domain, $f(x)$ represents a source or forcing term, and $\mathbf{A}: \Omega \times \mathbb{R}^d \to \mathbb{R}^d$ is a nonlinear flux vector field. 

From a broad physical perspective, equations of this structure transcend single disciplines:
\begin{itemize}
    \item \textbf{Solid mechanics:} The flux $\mathbf{A}(x, \nabla \Phi)$ models materials where stress or strain saturates near geometric stress concentrators, preventing unrealistic infinite stress singularities predicted by classical linear elasticity. Specifically, for anti-plane shear or Airy stress potentials $\Phi$, the constitutive flux takes the quasilinear form:
    \begin{equation}
    \mathbf{A}(x, \nabla \Phi) = \frac{\nabla \Phi}{2\mu(1 + \beta|\nabla \Phi|^\alpha)^{1/\alpha}},
    \end{equation}
    where $\mu > 0$ is the shear modulus, $\alpha > 1$ is the strain-limiting exponent controlling operator stiffness, and $\beta > 0$ is the nonlinearity modulus.
    \item \textbf{Fluid dynamics and porous media:} Similar quasilinear structures govern non-Darcy flows, such as the Darcy-Brinkman-Forchheimer models where quadratic drag terms introduce velocity-dependent resistance to fluid motion.
    \item \textbf{Electrohydrodynamics:} Charge-coupled fluid transport models incorporate nonlinear electric field interactions that modify momentum equations through complex source and damping dependencies.
\end{itemize}

The parameter $\alpha$ dictates how sharply the flux responds to growing gradients. As $\alpha$ increases, the operator transitions from a gentle nonlinearity to a stiff, highly localized response, creating severe gradients and challenging optimization landscapes.

\subsection{Constitutive restrictions}
The mathematical and physical richness of the quasilinear flux operator $\mathbf{A}(x, \nabla \Phi)$ stems from the theory of implicit constitutive relations, wherein the strain tensor is expressed as an implicit, nonlinear function of the stress tensor (or vice versa). In classical Cauchy elasticity, stress is a linear or moderately nonlinear explicit function of the strain gradient, which inherently leads to unbounded stress fields (singularities) at crack tips or sharp geometric notches. 

To resolve this unphysical behavior, strain-limiting theories posit that the maximum principal stretch or strain components are bounded, regardless of how large the applied loading becomes. For anti-plane shear deformations and potential-based formulations, this restriction manifests as a nonlinear algebraic coupling between the gradient norm $|\nabla \Phi|$ and the resulting stress field $\boldsymbol{\tau}$. The denominator term $(1 + \beta|\nabla \Phi|^\alpha)^{1/\alpha}$ acts as a self-adaptive attenuation mechanism: as the local gradient magnitude $|\nabla \Phi| \to \infty$, the effective compliance adjusts such that the stress saturates to a finite material limit. This structural property fundamentally alters the analytical character of the governing PDE, shifting it from uniformly elliptic (when $\alpha = 2$ and $\beta = 0$) to strongly degenerate or singular elliptic regimes depending on the choice of the nonlinearity exponent $\alpha > 0$.

\subsection{Monotonicity, coercivity, and generalized Sobolev spaces}
Rigorously establishing the existence and uniqueness of weak solutions for such strongly nonlinear problems requires moving beyond standard Hilbert space frameworks into generalized Lebesgue--Sobolev spaces. Because the growth conditions of the flux operator are dictated by the exponent $\alpha$, the natural energy space for the displacement or potential field is the reflexive Banach space $W^{1, p}(\Omega)$, where the integrability index $p$ is directly tied to the constitutive parameters.

The operator $\mathbf{A}: W^{1, p}(\Omega) \to [W^{1, p}(\Omega)]^*$ satisfies standard Leray--Lions type conditions:
\begin{enumerate}
    \item \textbf{Boundedness (Growth Condition):} There exist constants $C_1 > 0$ and $p > 1$ such that the flux satisfies polynomial growth bounds:
    \begin{equation}
    |\mathbf{A}(x, \xi)| \le C_1 \left( 1 + |\xi|^{p-1} \right) \quad \forall \xi \in \mathbb{R}^d.
    \end{equation}
    \item \textbf{Coercivity:} The operator exhibits strict coercivity to ensure boundedness of energy minimizing sequences:
    \begin{equation}
    \frac{\int_{\Omega} \mathbf{A}(x, \nabla \Phi) \cdot \nabla \Phi \, dx}{\|\Phi\|_{W^{1, p}(\Omega)}^{p}} \to \infty \quad \text{as } \|\Phi\|_{W^{1, p}(\Omega)} \to \infty.
    \end{equation}
    \item \textbf{Strict Monotonicity:} The operator satisfies the strict monotonicity condition, ensuring that the energy functional is strictly convex:
    \begin{equation}
    \int_{\Omega} \left( \mathbf{A}(x, \xi_1) - \mathbf{A}(x, \xi_2) \right) \cdot (\xi_1 - \xi_2) \, dx > 0 \quad \forall \xi_1 \neq \xi_2 \in \mathbb{R}^d.
    \end{equation}
\end{enumerate}
These structural properties guarantee that the associated calculus of variations minimization problem possesses a unique minimizer. However, the presence of sharp spatial transitions and high-gradient zones makes traditional grid-based discretizations prone to numerical dispersion or mesh locking, establishing the strong mathematical motivation for mesh-free, physics-informed neural network architectures that preserve global differentiability.

\section{Boundary value problem formulation and functional analysis}

To analyze the quasilinear boundary value problem rigorously, we augment the differential equation with homogeneous or inhomogeneous Dirichlet boundary conditions on the boundary $\partial\Omega$:
\begin{equation}
\begin{aligned}
-\nabla \cdot \left( \frac{\nabla \Phi}{2\mu(1 + \beta|\nabla \Phi|^\alpha)^{1/\alpha}} \right) &= f(x) && \text{in } \Omega, \\
\Phi &= g(x) && \text{on } \partial\Omega.
\end{aligned}
\end{equation}

\subsection{Variational setting and weak formulation}
We formulate the problem within appropriate Sobolev spaces. Let $H^1(\Omega)$ denote the standard Sobolev space of square-integrable functions possessing square-integrable weak derivatives, and let $H_0^1(\Omega)$ denote the subspace with zero trace on $\partial\Omega$. Multiplying the quasilinear divergence equation by an admissible test function $v \in H_0^1(\Omega)$ and applying the divergence theorem (integration by parts) yields the weak formulation:
\begin{equation}
\int_{\Omega} \frac{\nabla \Phi \cdot \nabla v}{2\mu(1 + \beta|\nabla \Phi|^\alpha)^{1/\alpha}} \, dx = \int_{\Omega} f v \, dx \quad \forall v \in H_0^1(\Omega).
\end{equation}

\subsection{Existence and uniqueness of solutions}
The mathematical analysis of quasilinear elliptic problems of this type relies heavily on the theory of monotone operators and the calculus of variations. 
\begin{theorem}[Existence and Uniqueness]
Under standard coercivity, hemicontinuity, and strict monotonicity conditions of the flux operator $\mathbf{A}(x, \xi)$ (satisfying polynomial growth conditions dictated by exponent $\alpha$), the weak formulation admits a unique solution $\Phi \in H^1(\Omega)$.
\end{theorem}
The proof typically utilizes the Minty-Browder theorem for monotone operators or direct minimization of the associated strictly convex energy functional in reflexive Banach spaces, ensuring that physical energy dissipation principles are mathematically well-posed. 

\subsection{The physical imperative for accurate numerical solutions}

Solving quasilinear equations accurately is not merely a numerical exercise; it is physically mandatory for several reasons:
\begin{enumerate}
    \item \textbf{Resolution of localization and boundary layers:} Quasilinear operators with saturation properties (such as strain-limiting or high Forchheimer numbers) frequently exhibit localized deformation zones or boundary layers. Inaccurate numerical approximations smear these gradients, masking the true physical mechanics of fracture initiation or fluid resistance.
    \item \textbf{Energy balance and constitutive admissibility:} Because the governing equations express local balance of linear momentum or mass coupled with nonlinear constitutive laws, pointwise residual minimization must faithfully reflect internal stress fields. An error in the gradient $\nabla \Phi$ directly distorts calculated stress tensors $\tau$, leading to unphysical predictions of material failure or flow rates.
    \item \textbf{Operator-solution decoupling:} As established in comparative studies, manufactured or physical source terms $f(x)$ depend explicitly on material parameters $(\alpha, \beta)$. Consequently, different parameter regimes solve distinct operators sharing identical exact solutions. High-accuracy solvers must isolate true solution error from operator-specific residual scaling.
\end{enumerate}

\section{Numerical method }
To approximate the solution of the quasilinear partial differential equation without traditional grid discretizations, we employ a mesh-free physics-informed framework powered by a Kolmogorov--Arnold Network. Unlike traditional multi-layer perceptrons that rely on fixed node-based activation functions and scalar connection weights, KANs position learnable one-dimensional functions directly on the network edges, offering superior representation efficiency and smoothness properties suitable for high-order differential operators.

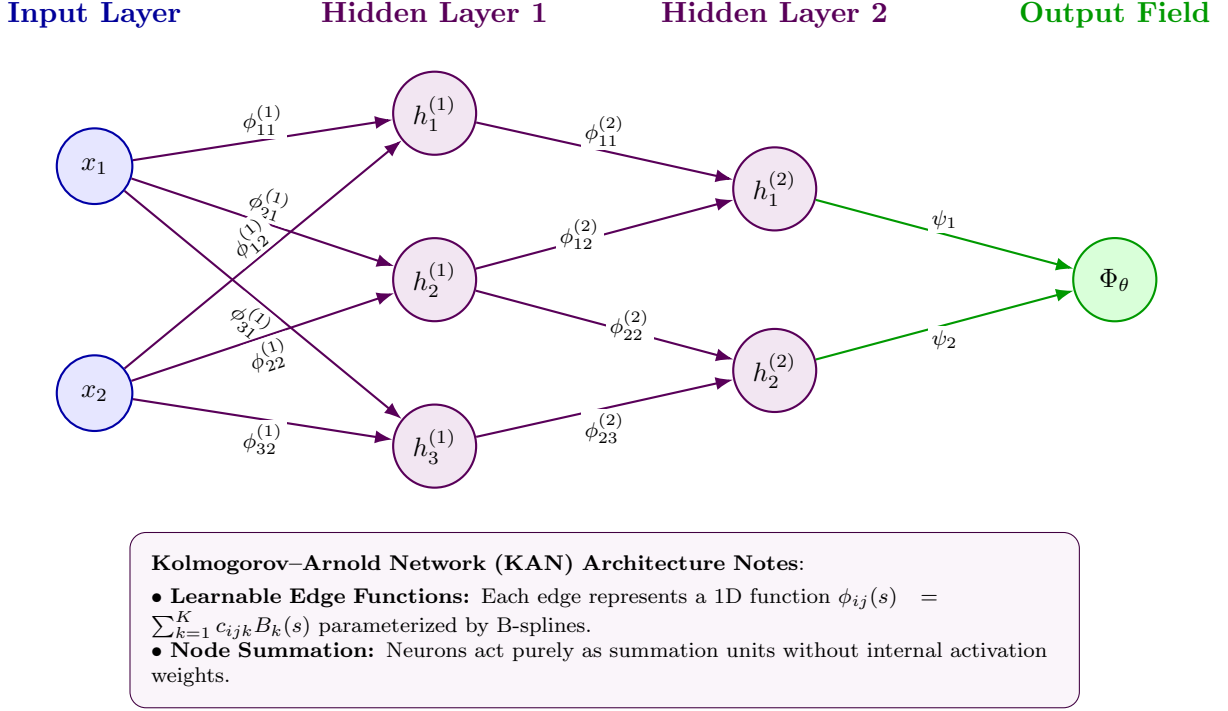
\begin{figure}[H]
\centering
\begin{tikzpicture}[
    >=Latex,
    font=\small,
    input node/.style={
        circle,
        draw=blue!65!black,
        fill=blue!10,
        minimum size=10mm,
        thick,
        align=center
    },
    hidden node/.style={
        circle,
        draw=violet!70!black,
        fill=violet!10,
        minimum size=10mm,
        thick,
        align=center
    },
    output node/.style={
        circle,
        draw=green!60!black,
        fill=green!15,
        minimum size=11mm,
        thick,
        align=center
    },
    edge label/.style={
        font=\scriptsize,
        fill=white,
        inner sep=1.5pt
    },
    box style/.style={
        rounded corners=6pt,
        draw=violet!50!black,
        fill=violet!4,
        align=left,
        inner sep=8pt,
        font=\scriptsize
    }
]

\node[font=\bfseries\color{blue!60!black}] at (0, 3.5) {Input Layer};
\node[font=\bfseries\color{violet!70!black}] at (4.5, 3.5) {Hidden Layer 1};
\node[font=\bfseries\color{violet!70!black}] at (9.0, 3.5) {Hidden Layer 2};
\node[font=\bfseries\color{green!60!black}] at (13.5, 3.5) {Output Field};


\node[input node] (x1) at (0,  1.5) {$x_1$};
\node[input node] (x2) at (0, -1.5) {$x_2$};

\node[hidden node] (h1_1) at (4.5,  2.2) {$h_1^{(1)}$};
\node[hidden node] (h1_2) at (4.5,  0.0) {$h_2^{(1)}$};
\node[hidden node] (h1_3) at (4.5, -2.2) {$h_3^{(1)}$};

\node[hidden node] (h2_1) at (9.0,  1.2) {$h_1^{(2)}$};
\node[hidden node] (h2_2) at (9.0, -1.2) {$h_2^{(2)}$};

\node[output node] (out) at (13.5, 0) {$\Phi_\theta$};


\draw[->, thick, draw=violet!70!black] (x1) -- (h1_1) node[midway, above, edge label] {$\phi_{11}^{(1)}$};
\draw[->, thick, draw=violet!70!black] (x1) -- (h1_2) node[midway, sloped, above, edge label] {$\phi_{21}^{(1)}$};
\draw[->, thick, draw=violet!70!black] (x1) -- (h1_3) node[midway, sloped, below, edge label] {$\phi_{31}^{(1)}$};

\draw[->, thick, draw=violet!70!black] (x2) -- (h1_1) node[midway, sloped, above, edge label] {$\phi_{12}^{(1)}$};
\draw[->, thick, draw=violet!70!black] (x2) -- (h1_2) node[midway, sloped, below, edge label] {$\phi_{22}^{(1)}$};
\draw[->, thick, draw=violet!70!black] (x2) -- (h1_3) node[midway, below, edge label] {$\phi_{32}^{(1)}$};

\draw[->, thick, draw=violet!70!black] (h1_1) -- (h2_1) node[midway, above, edge label] {$\phi_{11}^{(2)}$};
\draw[->, thick, draw=violet!70!black] (h1_2) -- (h2_1) node[midway, left, edge label] {$\phi_{12}^{(2)}$};
\draw[->, thick, draw=violet!70!black] (h1_2) -- (h2_2) node[midway, right, edge label] {$\phi_{22}^{(2)}$};
\draw[->, thick, draw=violet!70!black] (h1_3) -- (h2_2) node[midway, below, edge label] {$\phi_{23}^{(2)}$};

\draw[->, thick, draw=green!60!black] (h2_1) -- (out) node[midway, above, edge label] {$\psi_1$};
\draw[->, thick, draw=green!60!black] (h2_2) -- (out) node[midway, below, edge label] {$\psi_2$};

\node[box style, text width=12cm] at (6.75, -4.5) {
    \textbf{Kolmogorov--Arnold Network (KAN) Architecture Notes}:\\[1mm]
    $\bullet$ \textbf{Learnable Edge Functions:} Each edge represents a 1D function $\phi_{ij}(s) = \sum_{k=1}^{K} c_{ijk} B_k(s)$ parameterized by B-splines.\\
    $\bullet$ \textbf{Node Summation:} Neurons act purely as summation units without internal activation weights.
};

\end{tikzpicture}
\caption{Refined multi-layer architecture of the Kolmogorov--Arnold Network (KAN) for the strain-limiting potential field $\Phi_\theta(\mathbf{x})$, featuring explicit learnable univariate edge functions.}
\label{fig:refined_kan_architecture}
\end{figure}

\subsection{Kolmogorov--Arnold Network (KAN) formulation}
In conventional multi-layer perceptrons, nonlinear activations are applied locally at the nodes, and internal connections are scaled by static scalar weights. Conversely, KANs build upon the Kolmogorov--Arnold representation theorem by placing parameterized, learnable one-dimensional functions directly on the network edges. 

As illustrated in Figure~\ref{fig:refined_kan_architecture}, each connection $\phi_{oi}$ between nodes is modeled as a linear projection combined with a flexible expansion of radial basis functions (RBFs) or B-splines. Specifically, the univariate edge transformation is expressed as:
\begin{equation}
\phi_{oi}(u) = w_{oi}u + \sum_{k=1}^{K} c_{oik} \exp\left(-\left(\frac{u - \mu_k}{\sigma_k}\right)^2\right),
\end{equation}
where $\mu_k$ and $\sigma_k$ denote fixed centers and widths governing the basis functions, while $w_{oi}$ and $c_{oik}$ represent the trainable optimization parameters. The network topology comprises multiple hidden layers configured with a uniform width (e.g., 96 hidden units per layer, structured as $\mathbb{R}^2 \to \mathbb{R}^{96} \to \mathbb{R}^{96} \to \mathbb{R}^{96} \to \mathbb{R}$). Between consecutive hidden layers, smooth hyperbolic tangent ($\tanh$) activation functions are applied to ensure that high-order derivatives—crucial for computing the divergence operator $\nabla \cdot \tau$ via automatic differentiation—remain continuous and bounded across the entire domain.

\subsection{Exact boundary condition enforcement via hard-bc Ansatz}
To satisfy complex Dirichlet boundary conditions strictly and eliminate boundary residual errors during early training phases, the network output is formulated using a hard boundary condition ansatz:
\begin{equation}
\Phi(x,y) = E(x,y) + D(x,y) \cdot N_\theta(x,y),
\end{equation}
where $E(x,y)$ is an analytical function constructed to match the exact boundary data on $\partial\Omega$, $D(x,y)$ represents a boundary distance factor designed to vanish identically on all outer boundaries ($\left. D \right|_{\partial\Omega} = 0$), and $N_\theta(x,y)$ corresponds to the raw scalar output generated by the KAN model. This exact formulation guarantees that the boundary conditions are satisfied a priori, leaving the training procedure to focus exclusively on minimizing interior governing equation residuals.

\subsection{Physics-informed loss function and collocation strategy}
The interior collocation strategy relies on a mesh-free sampling approach, utilizing 2048 interior points resampled dynamically at each training epoch alongside 128 boundary collocation points per edge to maintain uniform spatial resolution. 

The governing physics-informed loss function is formulated by combining the interior residual of the quasilinear partial differential equation with a boundary constraint penalty. Let the nonlinear flux vector $\boldsymbol{\tau}(\mathbf{x})$ and the corresponding PDE residual $r(\mathbf{x}_i; \theta)$ at an interior collocation point $\mathbf{x}_i \in \Omega$ be defined as:
\begin{equation}
\boldsymbol{\tau}(\mathbf{x}_i) = \frac{\nabla \Phi_\theta(\mathbf{x}_i)}{2\mu \left(1 + \beta \|\nabla \Phi_\theta(\mathbf{x}_i)\|^\alpha \right)^{1/\alpha}},
\end{equation}
\begin{equation}
r(\mathbf{x}_i; \theta) = -\nabla \cdot \boldsymbol{\tau}(\mathbf{x}_i) - f(\mathbf{x}_i),
\end{equation}
where $\Phi_\theta(\mathbf{x}_i)$ represents the neural network approximation parameterized by weights and biases $\theta$, and $f(\mathbf{x}_i)$ is the given source term.

To ensure resilience against outliers and stiffness in the non-convex optimization landscape, the interior residual is penalized using a robust pseudo-Huber loss function:
\begin{equation}
\rho_{\delta}(z) = \delta^2 \left( \sqrt{1 + \left(\frac{z}{\delta}\right)^2} - 1 \right),
\end{equation}
where $\delta$ is the transition parameter controlling the crossover between $L_1$ and $L_2$-like penalty behaviors. 

Additionally, boundary conditions are enforced via a dedicated boundary loss term $\mathcal{L}_{\mathrm{bc}}$ evaluated over boundary collocation points $\mathbf{y}_j \in \partial\Omega$. The total composite loss function to be minimized is thus defined as:
\begin{equation}
\mathcal{L}(\theta) = \frac{\lambda_{\mathrm{pde}}}{N_u} \sum_{i=1}^{N_u} \rho_{\delta}\left(r(\mathbf{x}_i; \theta)\right) + \frac{\lambda_{\mathrm{bc}}}{N_b} \sum_{j=1}^{N_b} \left( \Phi_\theta(\mathbf{y}_j) - g(\mathbf{y}_j) \right)^2,
\end{equation}
where $N_u$ and $N_b$ denote the number of interior and boundary collocation points, respectively, while $\lambda_{\mathrm{pde}}$ and $\lambda_{\mathrm{bc}}$ serve as loss-balancing hyperparameters that weigh the relative importance of satisfying the governing differential equation versus the boundary data.

\subsection{Hybrid three-stage optimization protocol}
Due to the highly non-convex nature of physics-informed loss landscapes, training proceeds through a rigorous hybrid three-stage optimization schedule designed to balance rapid global descent with high-precision local refinement:
\begin{enumerate}
    \item \textbf{Coarse Adam Training (Stage 1):} Optimization begins with 2500 epochs using the Adam optimizer at a conservative learning rate of $1\times 10^{-3}$ to rapidly navigate toward the global basin of attraction.
    \item \textbf{Fine Adam Training (Stage 2):} The learning rate is reduced to $1\times 10^{-4}$ for an additional 3000 epochs to stabilize gradient updates and refine interior residual distributions.
    \item \textbf{L-BFGS Polish (Stage 3):} A final second-order quasi-Newton L-BFGS polish is executed for 200 epochs with strong Wolfe line searches to achieve high-precision convergence near the local minimum.
\end{enumerate}
To prevent collocation overfitting and ensure robust generalization, model selection is governed by evaluating performance against a fixed, held-out validation set comprising 2048 points, ensuring that the final restored checkpoint corresponds to the minimum validation loss rather than transient training fluctuations. 

\section{Training algorithm and KAN adaptivity mechanisms}

To solve the quasilinear strain-limiting partial differential equation without relying on static grids or manual discretization, the framework implements an advanced \textit{adaptive Kolmogorov-Arnold physics-informed neural networks} coupled with dynamic grid adaptation, basis expansion, and residual-based sampling strategies.

\subsection{Multi-stage training algorithm and optimization schedule}
The training process is formulated as a non-convex optimization problem solved through a robust, multi-stage protocol designed to balance global convergence stability with high-precision local refinement. 

\begin{enumerate}
    \item \textbf{Stage 1: Primary Adam optimization with curriculum ramping}\\
    Training initiates using the Adam optimizer with an exponential learning rate decay scheduler ($\gamma_{\mathrm{adam}}$). To prevent gradient pathologies caused by strong initial nonlinearities, a curriculum weighting schedule ramps the PDE residual weight $w_{\mathrm{pde}}(epoch)$ from a small initial value up to unity over a designated ramp period:
    \begin{equation}
    w_{\mathrm{pde}}(e) = \min\left(1.0, \max\left(0.0, \frac{e - e_{\mathrm{pre}}}{\epsilon_{\mathrm{ramp}}}\right)\right).
    \end{equation}
    
    \item \textbf{Stage 2: Finetuning and sharpening}\\
    The optimization transitions to a second phase utilizing either Adam or AdamW with weight decay ($\ell_2$ regularization) under a reduced learning rate ($\eta_{\mathrm{finetune}}$) to sharpen coefficient convergence and minimize interior residuals.
    
    \item \textbf{Stage 3: Second-order L-BFGS polish}\\
    To achieve high-precision convergence near the local minimum, a quasi-Newton L-BFGS optimizer is employed with strong Wolfe line searches. This stage incorporates curvature information (Hessian approximations via past gradient histories) and is conditionally triggered only if validation performance continues to improve.
\end{enumerate}

The composite loss function minimized across all stages aggregates interior governing equations, boundary constraints, and energy densities:
\begin{equation}
\mathcal{L}_{\mathrm{total}}(\theta) = \lambda_{\mathrm{pde}} \mathcal{L}_{\mathrm{pde}} + \lambda_{\mathrm{bc}} \mathcal{L}_{\mathrm{bc}} + \lambda_{\mathrm{energy}} \mathcal{L}_{\mathrm{energy}}.
\end{equation}

\subsection{KAN adaptivity and dynamic grid mechanisms}
Unlike standard multi-layer perceptrons with fixed internal structures, the KAN layers in the adaptive framework feature dynamic RBF/Gaussian basis centers and flexible grid topologies.

\begin{itemize}
\item \textbf{Dynamic grid relocation (center adaptation): }To concentrate basis functions where activation density is highest, internal grid centers $c_{ij}$ are periodically relocated toward the empirical distribution of layer inputs. For a given layer input sample set $S$, new centers are computed by blending uniform spacing with quantile-based adaptive spacing:
\begin{equation}
c_{\mathrm{new}} = \xi \, c_{\mathrm{uniform}} + (1 - \xi) \, c_{\mathrm{quantile}},
\end{equation}
where $\xi$ is a uniform mixing parameter (`$grid_mix_uniform$`). Kernel widths $\sigma_{ij}$ are automatically updated based on local inter-center distances scaled by a width factor.
\item \textbf{Progressive basis extension (grid growth): }To enhance network expressivity on-the-fly, the number of basis functions per layer ($\text{n\_basis}$) can grow dynamically from an initial size ($\text{n\_basis\_start}$) up to a maximum limit. When extension criteria are met:
\begin{itemize}
    \item The algorithm identifies spatial intervals with the largest gaps between existing RBF centers and inserts new midpoint centers.
    \item Coefficients for newly added basis functions are interpolated linearly from adjacent existing coefficients.
    \item Optimizer momentum buffers (such as Adam's first and second moment exponential moving averages, $m_t$ and $v_t$) are dynamically resized and mapped onto the enlarged parameter tensors to prevent training disruption.
\end{itemize}
\item \textbf{Residual-based adaptation: }To focus computational effort on regions with high mechanical stress gradients or unresolved boundary layers, the framework incorporates an adaptive collocation sampling scheme. Candidate interior points are sampled and evaluated for absolute PDE residual magnitude $|R_\theta(\mathbf{x})|$. New collocation points are sampled via probability weighting governed by a residual power law:
\begin{equation}
P(\mathbf{x}) \propto |R_\theta(\mathbf{x})|^{p} + \epsilon_{\mathrm{floor}},
\end{equation}
injecting high-error regions directly into the training set. Additionally, pointwise residual loss weighting dynamically scales individual collocation losses based on local error severity before backpropagation.
\end{itemize}

\section{Adaptive KAN-PINN solution methodology}

This investigation centers on the development of an advanced physics-informed Kolmogorov--Arnold Network (KAN-PINN) framework for approximating the solution to the multidimensional, quasilinear strain-limiting partial differential equation. A significant computational challenge in this context is the absence of a known closed-form analytical solution for complex geometries or general loading parameters, precluding traditional supervised learning paradigms that rely on exact target datasets. Our approach circumvents this limitation by embedding the governing quasilinear divergence operator and exact boundary constraints directly into a composite training objective.

Unlike conventional multilayer perceptrons that restrict nonlinear activations to static node operations and scale internal connections with scalar weights, KANs deploy learnable one-dimensional functions directly along the network edges. To maximize representation efficiency and capture sharp spatial variations without excessive model scaling, the proposed framework integrates several advanced algorithmic extensions: dynamic grid relocation (center adaptation), progressive basis extension (grid growth), residual-based adaptive collocation sampling, and exact boundary condition enforcement via a hard-BC ansatz.

A primary contribution of this work is a systematic characterization of the hyperparameter space and adaptivity configurations required to achieve robust, high-fidelity convergence. The sensitivity analysis was structured to rigorously evaluate the influence of the following core design elements:
\begin{enumerate}
    \item \textbf{Collocation point density and adaptive sampling:} The influence of uniform interior collocation density alongside residual-based adaptive point injection targeting high-error zones.
    \item \textbf{Basis function configuration and resolution:} The impact of the number of basis functions ($K$) and kernel width scaling on the expressivity of the learnable edge functions.
    \item \textbf{Network topology (depth and width):} The representational capacity governed by the number of hidden layers and layer width (neurons per layer).
    \item \textbf{Hard boundary condition enforcement:} The efficacy of exact Dirichlet satisfaction versus soft penalty weighting.
    \item \textbf{Multi-stage hybrid optimization schedule:} The convergence dynamics and stability provided by transitioning from stochastic gradient methods (Adam/AdamW) to second-order quasi-Newton polishing (L-BFGS).
\end{enumerate}

\subsection{Adaptive grid relocation and basis function dynamics}
A defining feature of the proposed RBF-KAN architecture is its reliance on grid-dependent Gaussian radial basis functions (RBFs) placed on network edges. Unlike fixed-grid methods, the internal RBF centers $\mu_k$ evolve dynamically in response to the feature distributions of layer inputs. 

During training, periodic grid adaptation updates the center locations by blending uniform spacing with quantile-based empirical distributions:
\begin{equation}
\mu_k^{\text{new}} = \xi \, \mu_k^{\text{uniform}} + (1 - \xi) \, \mu_k^{\text{quantile}},
\end{equation}
where $\xi$ is the uniform mixing parameter (`$grid_mix_uniform$`). Kernel widths $\sigma_k$ are concurrently scaled relative to local inter-center distances, ensuring high resolution where gradients are steep and broad smoothing where fields are smooth. Furthermore, progressive basis extension allows the network to start with a compact basis size ($K_{\text{start}}$) and grow incrementally up to a target resolution ($K_{\text{max}}$) by inserting midpoint centers in regions of maximum spacing. To maintain training momentum during grid updates, optimizer moment buffers (e.g., Adam exponential moving averages) are dynamically mapped onto the expanded coefficient tensors.

\subsection{Residual-based adaptive collocation sampling}
To prevent under-resolution in localized deformation zones, uniform interior sampling is augmented with residual-based adaptive collocation. At designated training intervals, candidate points are sampled across the domain $\Omega$ and evaluated for local PDE residual magnitude $|R_\theta(\mathbf{x})|$. New collocation points are injected into the training batch via probability-weighted sampling governed by a residual power law:
\begin{equation}
P(\mathbf{x}) = \frac{|R_\theta(\mathbf{x})|^p + \epsilon_{\text{floor}}}{\sum_{j} \left(|R_\theta(\mathbf{x}_j)|^p + \epsilon_{\text{floor}}\right)},
\end{equation}
thereby directing computational gradient tracking toward regions experiencing severe constitutive nonlinearity or elevated stress concentrations.

\subsection{Hybrid three-stage optimization protocol}
Due to the non-convex nature of quasilinear PINN loss landscapes, training relies on a robust three-stage hybrid optimizer schedule:
\begin{itemize}
    \item \textbf{Stage 1 (coarse adam training):} Optimization begins with the Adam optimizer at a higher learning rate ($\eta = 1\times 10^{-3}$) paired with a PDE residual curriculum ramp ($w_{\text{pde}} \to 1$) to rapidly navigate toward the global basin of attraction[cite: 1, 2, 3].
    \item \textbf{Stage 2 (finetuning and regularization):} The learning rate is decayed (e.g., $\eta = 1\times 10^{-4}$), and optional weight decay or AdamW regularization is applied to suppress high-frequency noise and refine interior residual distributions.
    \item \textbf{Stage 3 (L-BFGS Quasi-Newton polish):} A second-order L-BFGS optimizer with strong Wolfe line searches is executed to achieve high-precision convergence near the local minimum, conditioned on continuous validation improvement.
\end{itemize}

Model selection is strictly governed by evaluating performance against an independent validation set, ensuring the final restored checkpoint corresponds to optimal generalization rather than transient training fluctuations.

\section{Numerical solutions and discussion}

To evaluate the efficacy, adaptability, and robustness of the proposed AKAPINN framework, controlled numerical experiments were executed across three distinct material parameter configurations: 
\begin{itemize}
\item Run 1 with parameter values $\alpha = 1.5, \beta = 0.4$
\item Run 2  with parameter values $\alpha = 2.4, \beta = 0.3$
\item Run 3 with parameter values $\alpha = 3.1, \beta = 0.2$.
\end{itemize}
These parameter values are verified through a rigorous machine learning approach developed in \cite{pati2026sle}, and have been earlier utilized in several numerical studies  \cite{vasilyeva2024generalized}. Each configuration solves the quasilinear strain-limiting partial differential equation on the unit square domain with the manufactured verification solution $\Phi_{\text{exact}}(x,y) = (\frac{\mu}{2})y^2$ \cite{YoonMallikarjunaiah2022}. All computational training runs, dynamic grid adaptations, streaming chunk-wise backpropagation steps, and field diagnostic exports were executed locally on an Apple MacBook Pro powered by an \textit{Apple Silicon processor}, leveraging its high-bandwidth \textit{Unified Memory Architecture} and the \textit{PyTorch Metal Performance Shaders} backend to accelerate tensor operations and automatic differentiation.

\subsection{Computational workflow and algorithm implementation}
The end-to-end execution of the mesh-free solver follows the systematic procedure outlined in Algorithm~\ref{alg:akapinn}. Training initiates with device configuration and model instantiation, leveraging an initial basis size ($n_{\text{basis\_start}}$) coupled with an exact hard boundary condition ansatz to eliminate boundary residual contamination a priori. The optimization pipeline is structured around a hybrid three-stage protocol: initial coarse descent via the Adam optimizer, fine-tuning under reduced learning rates, and a final second-order quasi-Newton L-BFGS polish utilizing strong Wolfe line searches. Throughout training, dynamic RBF center relocation and progressive basis extension update network expressivity in high-gradient zones, while chunk-wise gradient accumulation manages memory efficiency on local hardware.

\begin{algorithm}[H]
\caption{Adaptive Grid RBF-KAN PINN Training Pipeline}
\label{alg:akapinn}
\begin{algorithmic}[1]
\Require Material parameters ($\mu, \beta, \alpha$), Geometry ($x_{\min}, x_{\max}, y_{\min}, y_{\max}$), Training hyperparameters ($\text{trn}$).
\Ensure Optimal model weights, validation logs, and field diagnostics.

\State Initialize Device (CPU/GPU) and set random seeds ($\text{trn.seed}$).
\State Instantiate Model:
\State \quad a. Create \texttt{AdaptiveGridRBFKANPINN} with hidden dimension and initial basis size (\texttt{n\_basis\_start}).
\State \quad b. Configure boundary ansatz parameters.
\State Build Training Stages \& Epoch Bounds:
\State \quad a. Stage 1: Adam optimizer (\texttt{learning\_rate}, \texttt{lr\_gamma\_adam}).
\State \quad b. Stage 2: Adam/AdamW optimizer (\texttt{finetune\_lr}, \texttt{weight\_decay}).
\State \quad c. Stage 3: L-BFGS optimizer (\texttt{lbfgs\_lr}, \texttt{history\_size}).
\State \textbf{If} $\text{trn.fixed\_collocation}$ is True:
\State \quad Sample fixed interior collocation points and boundary points.
\State \textbf{For} each training stage in [Stage 1, Stage 2, Stage 3]:
\State \quad \textbf{For} each $\text{epoch}$ from $\text{stage\_start}$ to $\text{stage\_end}$:
\State \quad \quad Update $\alpha$ via homotopy continuation schedule.
\State \quad \quad Sample or retrieve interior and boundary collocation points.
\State \quad \quad \textbf{If} grid adaptation/extension is enabled and epoch criteria met:
\State \quad \quad \quad Adapt RBF centers to empirical input distributions.
\State \quad \quad \quad Extend basis size up to \texttt{n\_basis} limit and update optimizer state buffers.
\State \quad \quad Compute Forward Pass \& Losses:
\State \quad \quad \quad Calculate Boundary Loss ($\mathcal{L}_{\text{bc}}$).
\State \quad \quad \quad Calculate Energy Loss ($\mathcal{L}_{\text{energy}}$) if enabled.
\State \quad \quad \quad Calculate PDE Residual Loss ($\mathcal{L}_{\text{pde}}$) via streaming chunk backward passes with optional point-wise residual weighting.
\State \quad \quad Perform Backpropagation \& Optimizer Step:
\State \quad \quad \quad Accumulate total loss: $\mathcal{L}_{\text{total}} = \lambda_{\text{bc}}\mathcal{L}_{\text{bc}} + \lambda_{\text{pde}}\mathcal{L}_{\text{pde}} + \lambda_{\text{energy}}\mathcal{L}_{\text{energy}}$.
\State \quad \quad \quad Clip gradients by norm (\texttt{max\_grad\_norm}).
\State \quad \quad \quad Step optimizer and learning rate scheduler.
\State \quad \quad Evaluate Validation Metrics:
\State \quad \quad \quad Compute validation loss at regular intervals (\texttt{validation\_every}).
\State \quad \quad \quad Track best model checkpoint based on validation performance (\texttt{min\_improve}).
\State \quad \quad Check Early Stopping criteria (\texttt{patience}).
\State Restore best-performing model checkpoint weights.
\State Execute Final Verification, Export Field Diagnostics (\texttt{.npz}), and Save Loss History / Plots.
\end{algorithmic}
\end{algorithm}

\subsection{Field approximation and boundary enforcement}
\begin{figure}[H]
    \centering
    \begin{minipage}{0.48\textwidth}
        \centering
        \includegraphics[width=\linewidth]{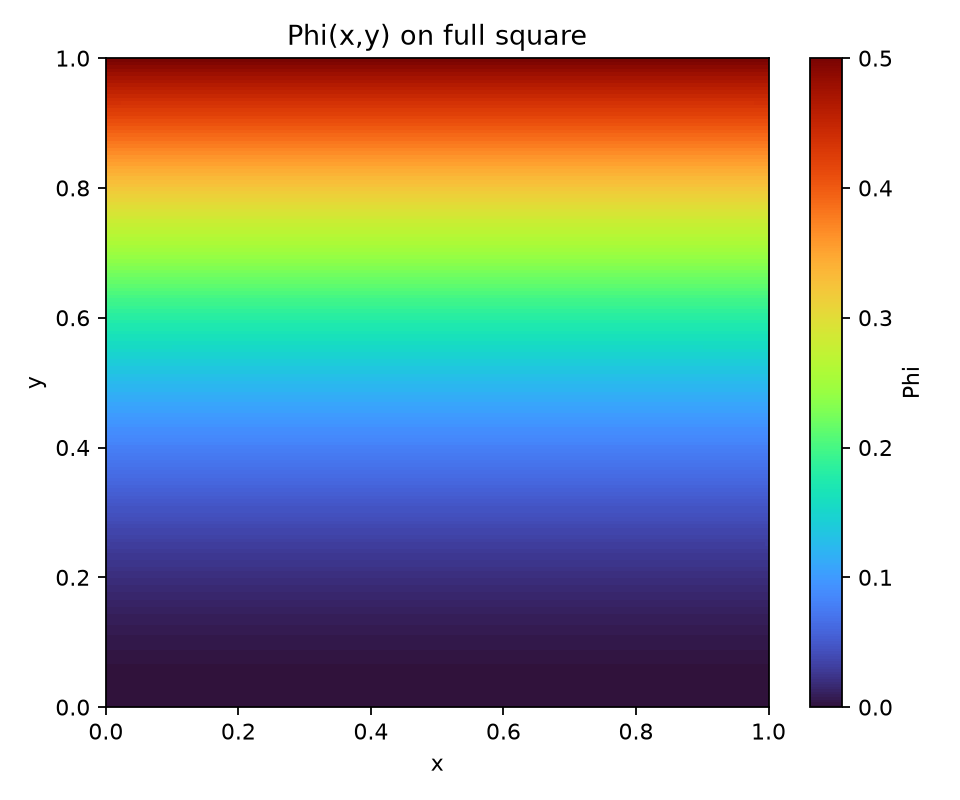}
        \caption{Run 1 ($\alpha=1.5, \beta=0.4$): Reconstructed potential field $\Phi(x,y)$ on the full square domain.}
        \label{fig:run1_phi}
    \end{minipage}\hfill
    \begin{minipage}{0.48\textwidth}
        \centering
        \includegraphics[width=\linewidth]{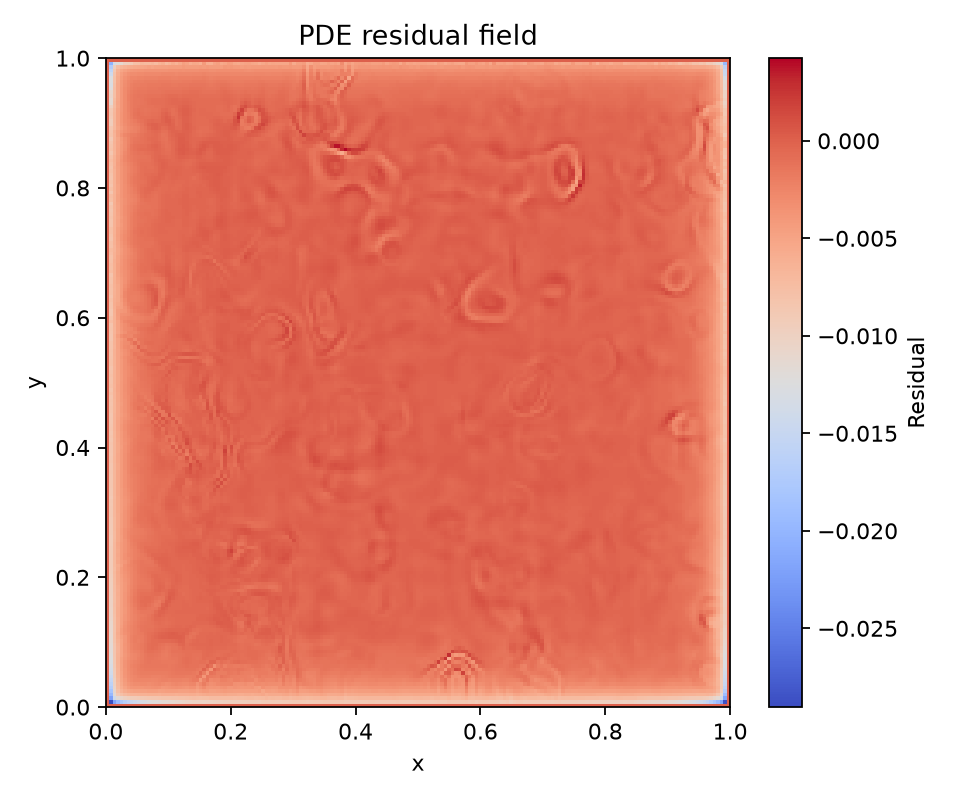}
        \caption{Run 1: PDE residual field illustrating spatial error distribution and operator scaling.}
        \label{fig:run1_res}
    \end{minipage}
\end{figure}

The reconstructed potential fields $\Phi(x,y)$ across all three runs exhibit excellent qualitative and quantitative agreement with the analytical solution, maintaining sub-$0.02\%$ relative $L_2$ error and identically zero boundary error. Figures~\ref{fig:run1_phi}, \ref{fig:run2_phi}, and \ref{fig:run3_phi} display the spatial distribution of $\Phi(x,y)$ for Runs 1, 2, and 3, respectively. Because the hard boundary condition ansatz strictly enforces Dirichlet data a priori, boundary discrepancies are entirely eliminated, allowing the hybrid three-stage optimization protocol (Adam to L-BFGS) to focus exclusively on minimizing interior governing equation residuals.

\subsection{Residual distributions and operator stiffness}
\begin{figure}[H]
    \centering
    \begin{minipage}{0.48\textwidth}
        \centering
        \includegraphics[width=\linewidth]{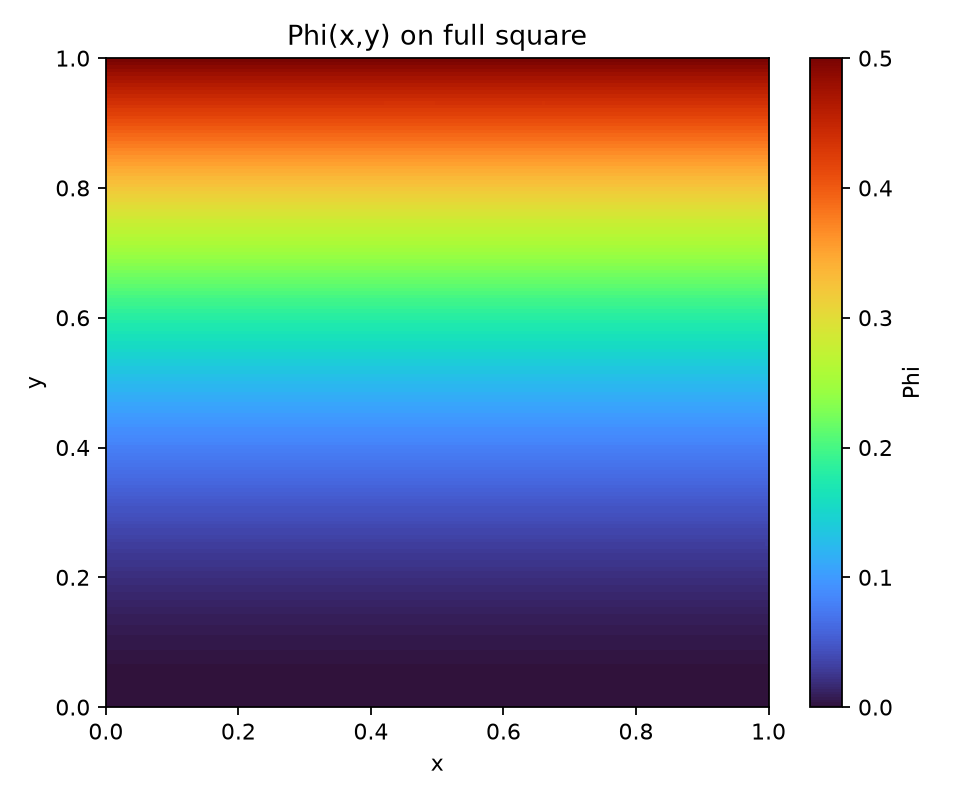}
        \caption{Run 2 ($\alpha=2.4, \beta=0.3$): Reconstructed potential field $\Phi(x,y)$.}
        \label{fig:run2_phi}
    \end{minipage}\hfill
    \begin{minipage}{0.48\textwidth}
        \centering
        \includegraphics[width=\linewidth]{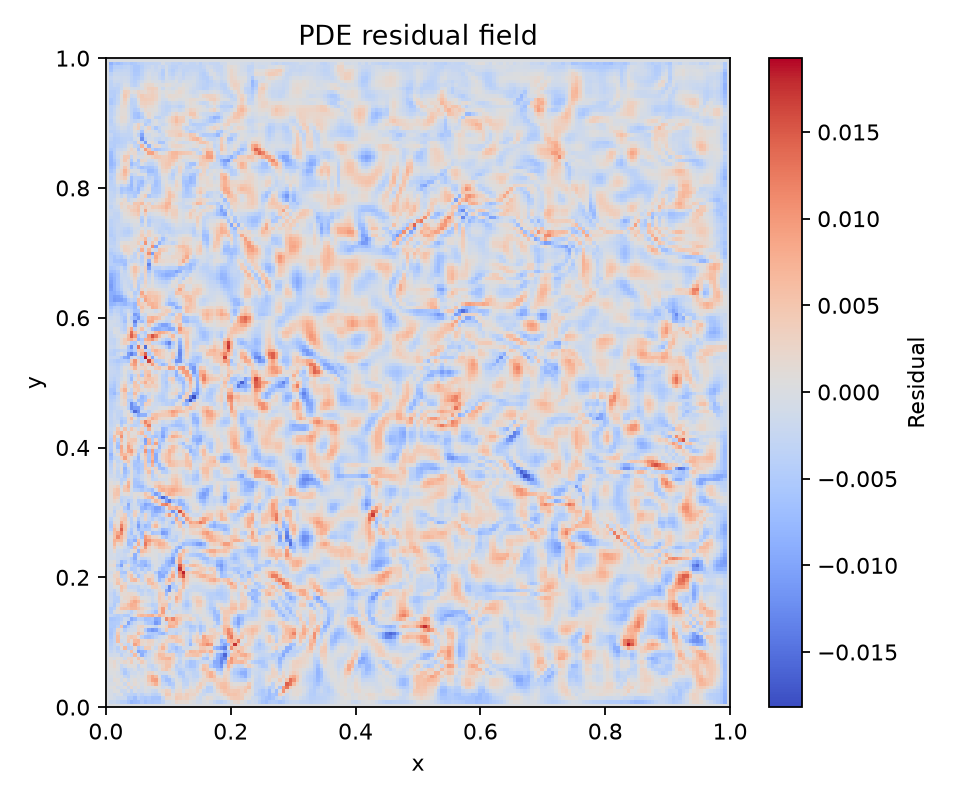}
        \caption{Run 2: Corresponding PDE residual field across the collocation grid.}
        \label{fig:run2_res}
    \end{minipage}
\end{figure}

\begin{figure}[H]
    \centering
    \begin{minipage}{0.48\textwidth}
        \centering
        \includegraphics[width=\linewidth]{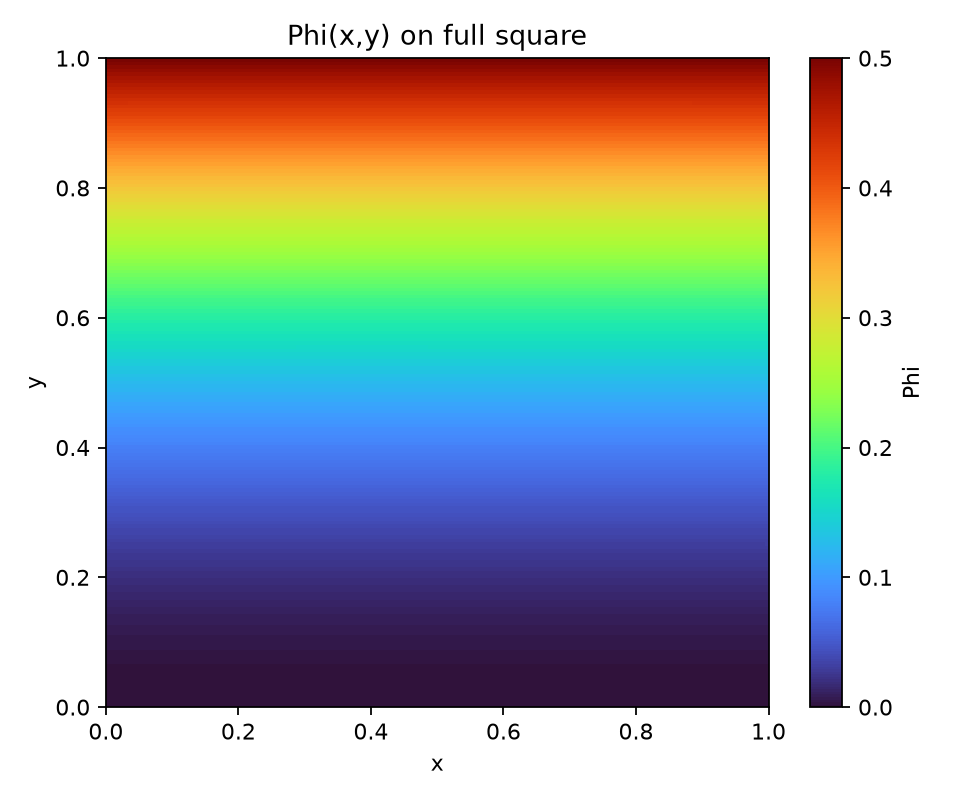}
        \caption{Run 3 ($\alpha=3.1, \beta=0.2$): Reconstructed potential field $\Phi(x,y)$.}
        \label{fig:run3_phi}
    \end{minipage}\hfill
    \begin{minipage}{0.48\textwidth}
        \centering
        \includegraphics[width=\linewidth]{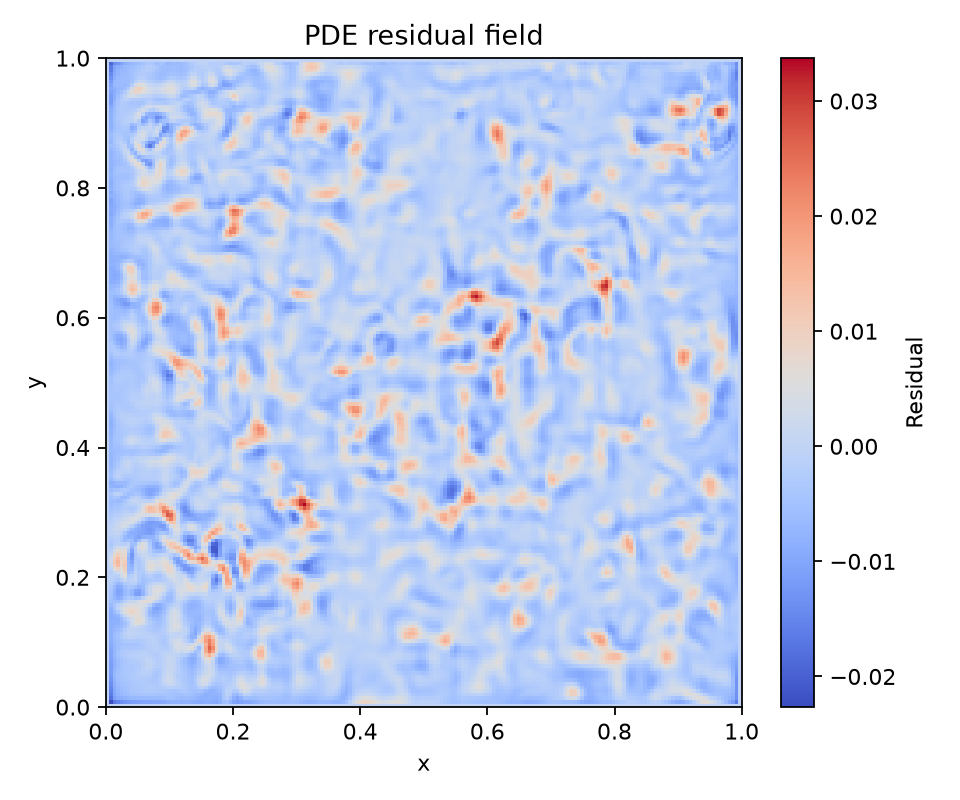}
        \caption{Run 3: Corresponding PDE residual field highlighting localized error structures.}
        \label{fig:run3_res}
    \end{minipage}
\end{figure}

Figures~\ref{fig:run1_res}, \ref{fig:run2_res}, and \ref{fig:run3_res} illustrate the corresponding interior PDE residual fields $r(x,y;\theta)$ for each parameter set. The spatial distribution of residuals highlights the influence of the strain-limiting exponent $\alpha$ and nonlinearity modulus $\beta$ on operator stiffness. In Run 1 ($\alpha=1.5$), the milder nonlinearity produces smoother flux responses, resulting in a lower final training loss ($1.15 \times 10^{-7}$). Conversely, higher exponents in Run 2 ($\alpha=2.4$) and Run 3 ($\alpha=3.1$) introduce sharper gradient sensitivities, shifting the terminal training loss values to $1.30 \times 10^{-6}$ and $4.59 \times 10^{-7}$, respectively.

\subsection{Metric inversion and generalization analysis}
A critical finding from this comparative study is the observation of a fundamental \emph{metric inversion} between training loss and true solution accuracy. Although Run 1 achieves the lowest training loss due to operator-specific residual scaling and lower constitutive stiffness, evaluation via the operator-independent relative $L_2$ error reveals that intermediate parameter regimes deliver superior generalization. Specifically, Run 2 achieves the highest accuracy with a relative $L_2$ error of $3.54 \times 10^{-5}$, followed by Run 3 ($7.78 \times 10^{-5}$) and Run 1 ($1.12 \times 10^{-4}$). Furthermore, Run 2 exhibits the tightest validation-to-training generalization gap ($1.48\times$), confirming that relative $L_2$ error—rather than cross-operator training loss—serves as the reliable model selection metric for quasilinear physics-informed neural networks.

\subsection{Training loss dynamics and convergence history}

\begin{figure}[H]
    \centering
    \begin{subfigure}[b]{0.32\linewidth}
        \centering
        \includegraphics[width=\linewidth]{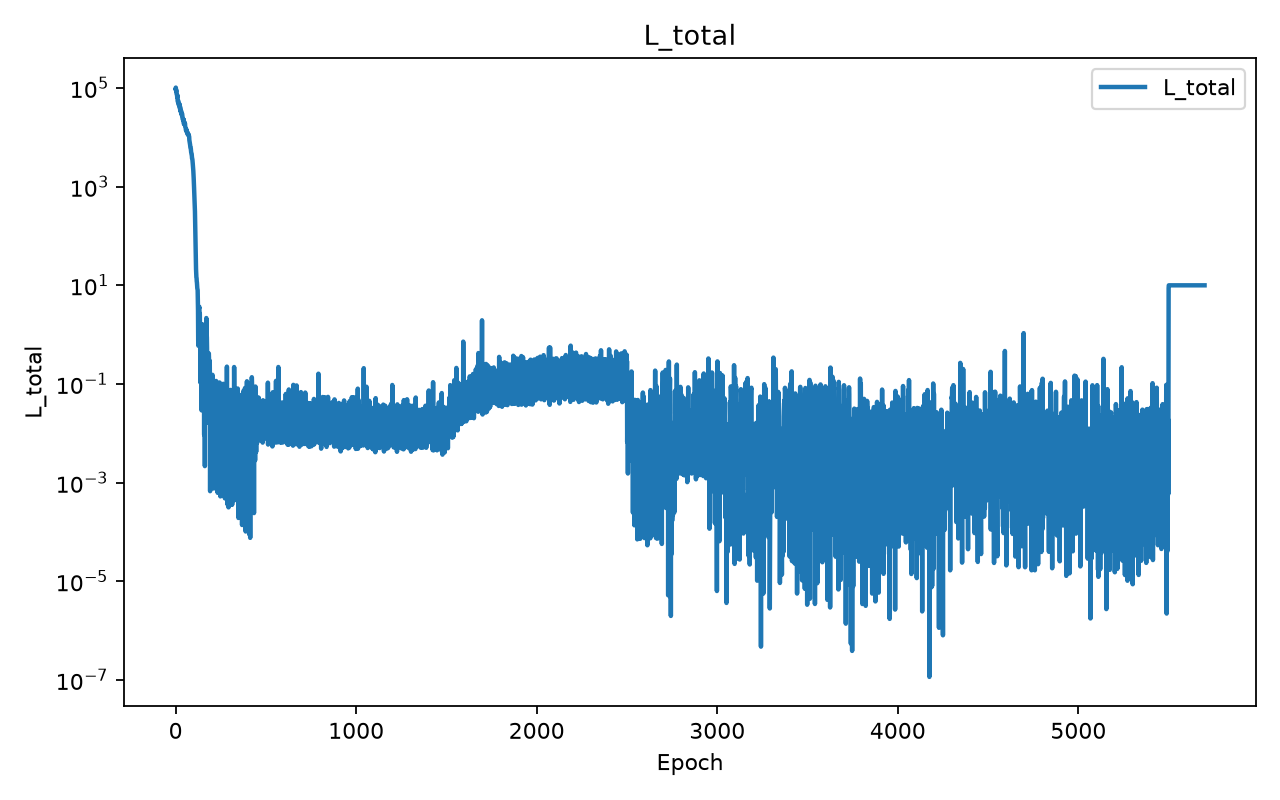}
        \caption{Run 1 ($\alpha = 1.5, \beta = 0.4$)}
        \label{fig:run1_loss}
    \end{subfigure}%
    \hfill
    \begin{subfigure}[b]{0.32\linewidth}
        \centering
        \includegraphics[width=\linewidth]{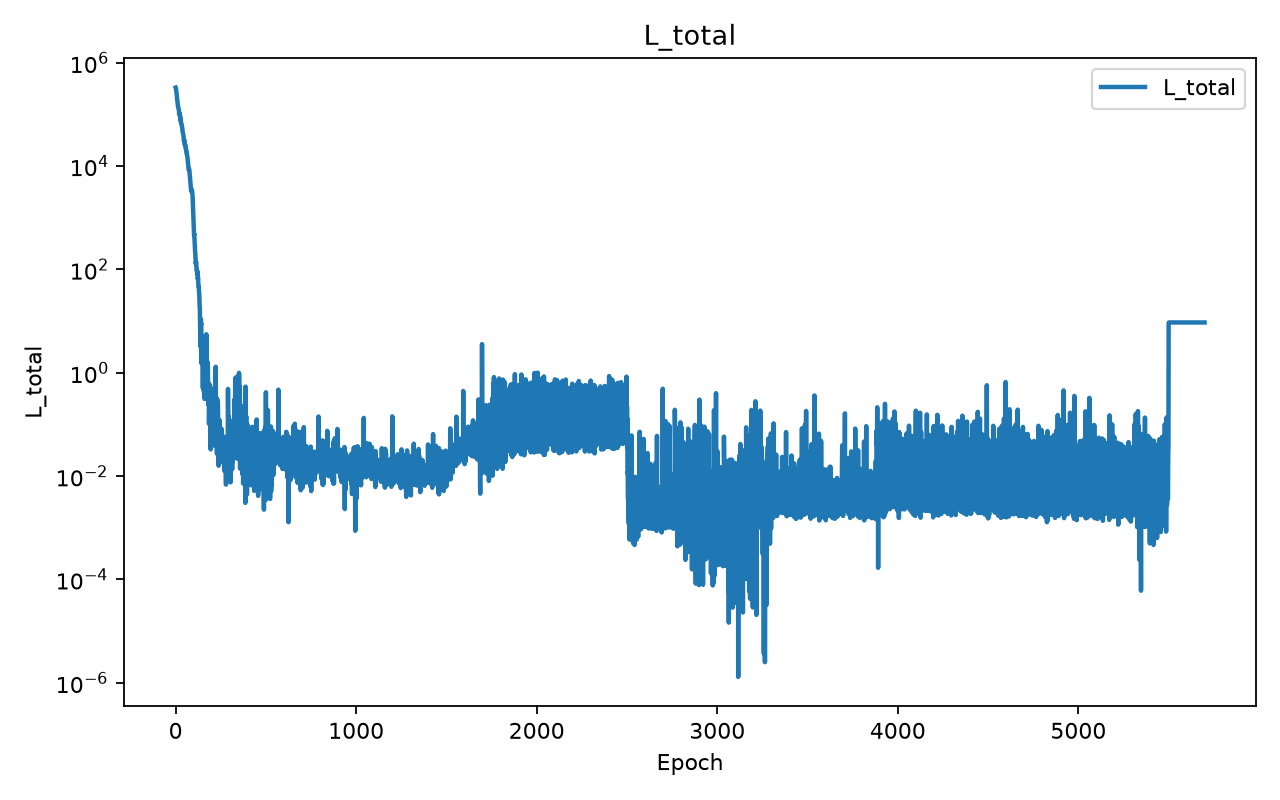}
        \caption{Run 2 ($\alpha = 2.4, \beta = 0.3$)}
        \label{fig:run2_loss}
    \end{subfigure}%
    \hfill
    \begin{subfigure}[b]{0.32\linewidth}
        \centering
        \includegraphics[width=\linewidth]{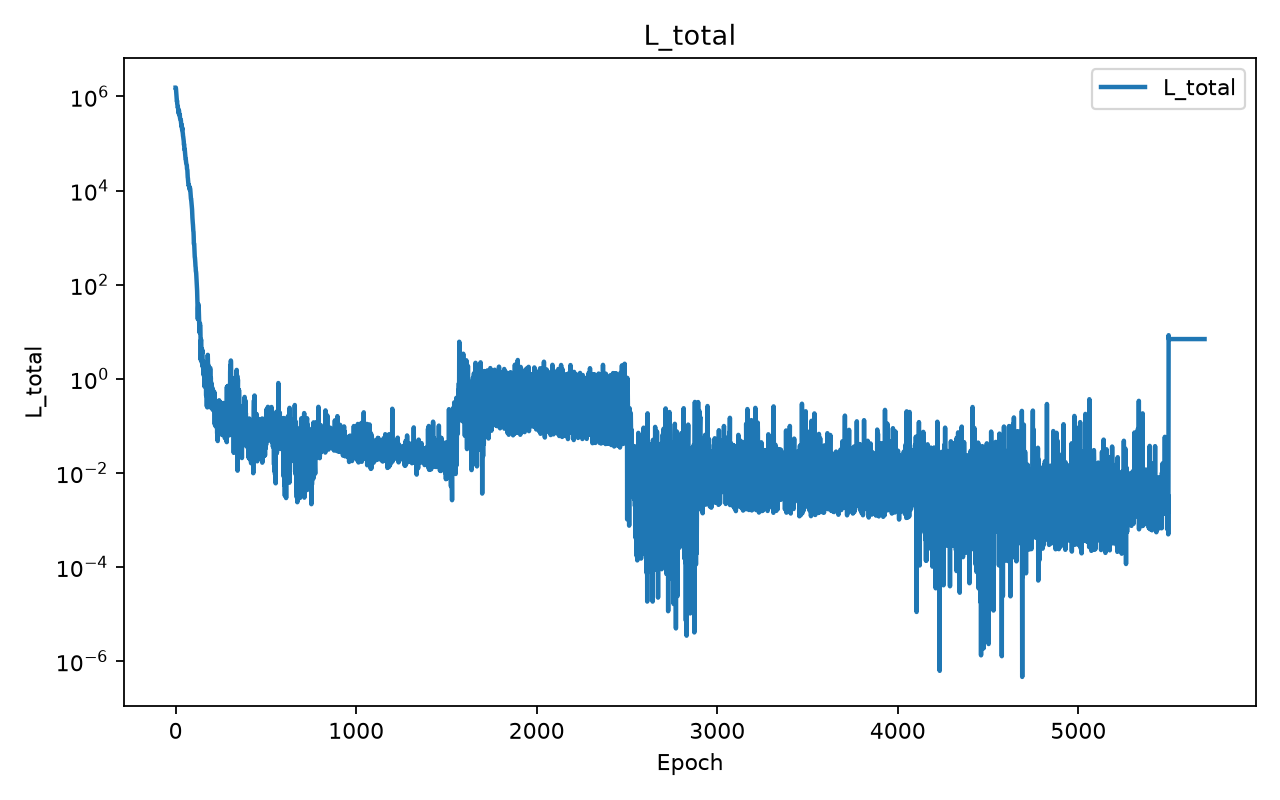}
        \caption{Run 3 ($\alpha = 3.1, \beta = 0.2$)}
        \label{fig:run3_loss}
    \end{subfigure}
    
    \caption{Total training loss history $\mathcal{L}_{\text{total}}$ across the hybrid three-stage optimization schedule for the three material parameter configurations .}
    \label{fig:all_loss_histories}
\end{figure}

The evolution of the total composite loss $\mathcal{L}_{\text{total}}$ across the three-stage optimization protocol highlights the convergence characteristics and robustness of the AKAPINN framework when executed locally on the Apple MacBook Pro. Figures~\ref{fig:run1_loss}, \ref{fig:run2_loss}, and \ref{fig:run3_loss} illustrate the loss histories for Run 1 ($\alpha = 1.5, \beta = 0.4$), Run 2 ($\alpha = 2.4, \beta = 0.3$), and Run 3 ($\alpha = 3.1, \beta = 0.2$), respectively. 

Across all three configurations, training initiates with a rapid, steep descent during the initial coarse Adam epochs, reducing the composite loss from initial magnitudes exceeding $10^5$ down to approximately $10^{-2}$ to $10^{-3}$ . An intermediate stabilization phase follows as the learning rate is refined and curriculum weights adjust. Upon transitioning to the second-order L-BFGS quasi-Newton polish phase (near epoch 2500), a sharp drop in residual magnitude occurs, allowing the optimizer to navigate complex non-convex basins and achieve terminal training loss values of $1.15 \times 10^{-7}$ for Run 1, $1.30 \times 10^{-6}$ for Run 2, and $4.59 \times 10^{-7}$ for Run 3. As established by the operator scaling, Run 1 achieves the lowest training loss due to its milder nonlinearity and lower stiffness, whereas higher values of the strain-limiting exponent $\alpha$ in Runs 2 and 3 introduce elevated gradient sensitivities in the loss landscape.

\begin{figure}[H]
    \centering
    \includegraphics[width=\linewidth]{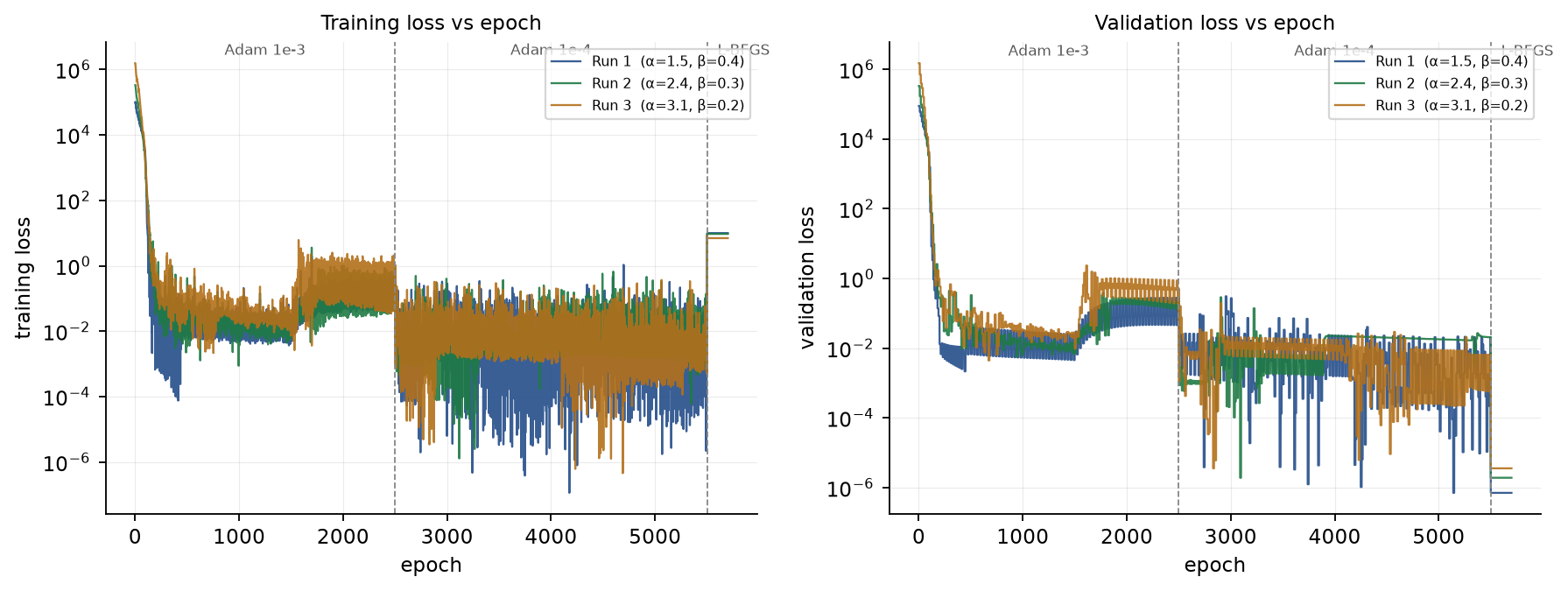}
    \caption{Overlaid training loss (left) and validation loss (right) trajectories versus epoch across the three material parameter configurations, highlighting distinct stage transitions (Adam $1\text{e-}3$, Adam $1\text{e-}4$, and L-BFGS) and the resulting metric inversion .}
    \label{fig:loss_trajectories}
\end{figure}

Overlaid loss trajectories across the three material parameter configurations—Run 1 ($\alpha=1.5, \beta=0.4$), Run 2 ($\alpha=2.4, \beta=0.3$), and Run 3 ($\alpha=3.1, \beta=0.2$)—demonstrate a unified convergence signature throughout the hybrid optimization protocol. As illustrated in Figure~\ref{fig:loss_trajectories}, all three runs share identical descent dynamics: a steep initial drop during the coarse Adam phase ($1\times10^{-3}$), followed by continued stabilization in the fine Adam phase ($1\times10^{-4}$), and a characteristic high-frequency fluctuation tail culminating in the L-BFGS quasi-Newton polish. The trajectories diverge primarily near the convergence threshold, revealing a critical metric inversion: although Run 1 achieves the lowest training loss due to smoother operator-specific residual scaling, it does not yield the optimal validation-set performance. Instead, intermediate regimes such as Run 2 exhibit superior generalization and lower true relative $L_2$ error, underscoring that training loss minimization alone is insufficient for model selection in quasilinear physics-informed neural networks. 

\subsection{The metric inversion phenomenon and operator stiffness}

A central insight revealed by the comparative analysis is the presence of a fundamental \emph{metric inversion} between training loss minimization and true solution accuracy. Placing the training loss distribution directly alongside the relative $L_2$ error underscores this contrast: the parameter configuration that minimizes the training loss is, paradoxically, the least accurate against the exact solution. 

\begin{figure}[H]
    \centering
    \begin{minipage}{0.48\textwidth}
        \centering
        \includegraphics[width=\linewidth]{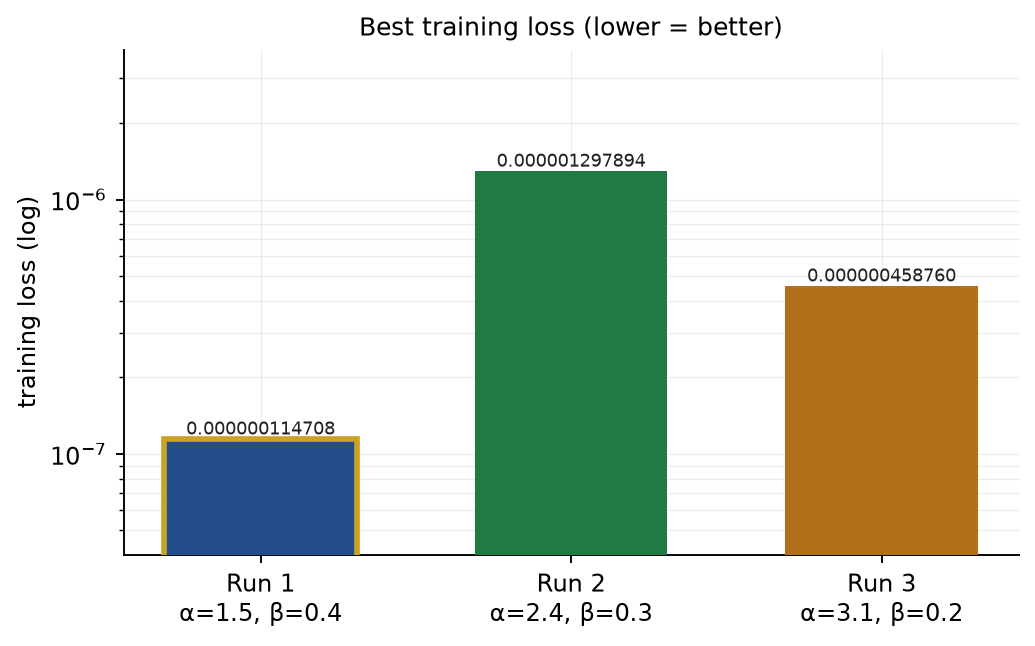}
        \caption{Best training loss across the three parameter configurations, showing Run 1 achieving the lowest value.}
        \label{fig:bar_training_loss}
    \end{minipage}\hfill
    \begin{minipage}{0.48\textwidth}
        \centering
        \includegraphics[width=\linewidth]{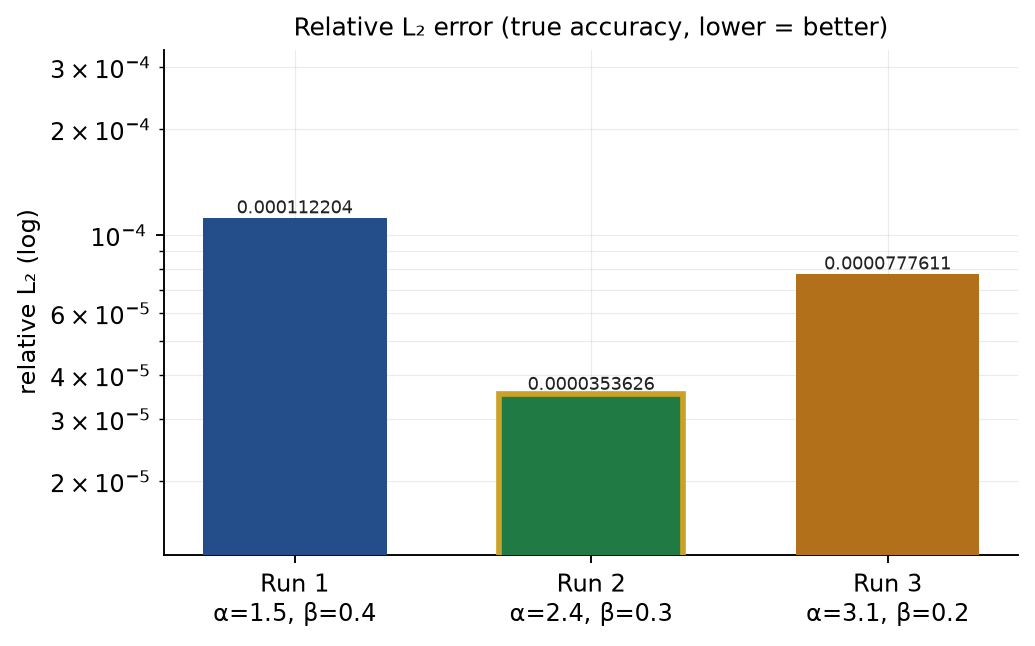}
        \caption{Relative $L_2$ error against the exact solution, where the performance bars flip and Run 2 achieves peak accuracy.}
        \label{fig:bar_relative_l2}
    \end{minipage}
\end{figure}

As illustrated in Figures~\ref{fig:bar_training_loss} and \ref{fig:bar_relative_l2}, Run 1 ($\alpha = 1.5, \beta = 0.4$) achieves the lowest training loss ($1.15 \times 10^{-7}$) due to its smoother operator-specific residual scaling and lower constitutive stiffness. However, when evaluated against the true manufactured solution via the operator-independent relative $L_2$ error, Run 1 yields the highest error ($1.12 \times 10^{-4}$). Conversely, Run 2 ($\alpha = 2.4, \beta = 0.3$) exhibits a higher training loss ($1.30 \times 10^{-6}$) yet achieves superior solution accuracy, attaining a peak relative $L_2$ error of $3.54 \times 10^{-5}$. 

\begin{figure}[H]
    \centering
    \includegraphics[width=0.85\linewidth]{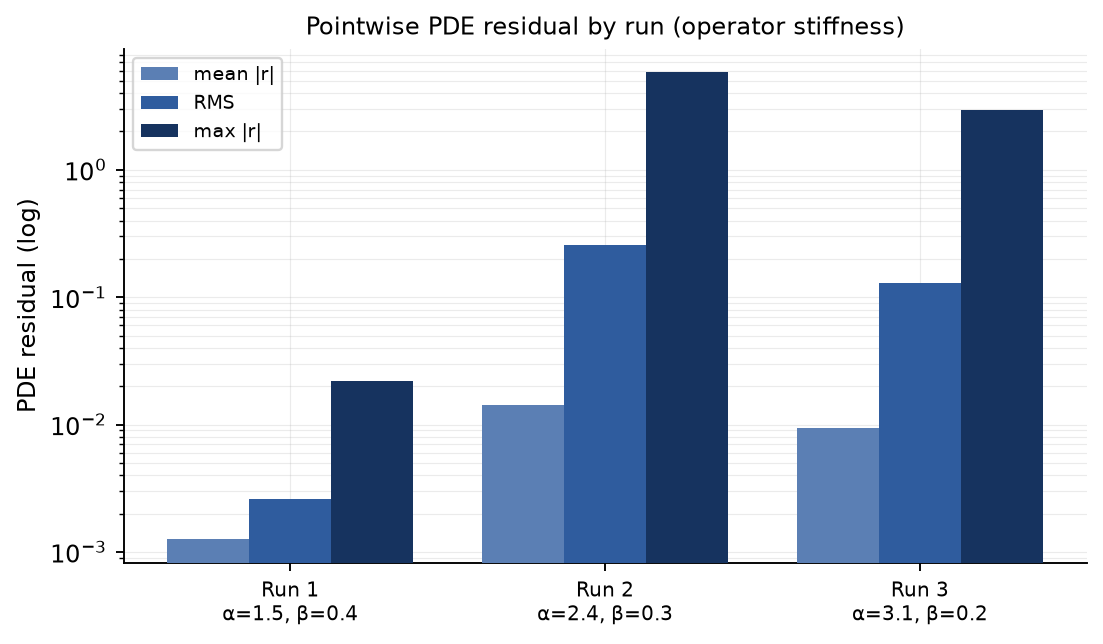}
    \caption{Pointwise PDE residual metrics (mean $|r|$, RMS, and max $|r|$) on a logarithmic scale, illustrating how Run 2's operator produces significantly larger residual magnitudes due to increased stiffness.}
    \label{fig:residual_bars}
\end{figure}

This reversal is directly explained by examining the pointwise PDE residual statistics shown in Figure~\ref{fig:residual_bars}. The operator governing Run 2 produces residual spikes (mean, RMS, and maximum absolute values) orders of magnitude larger than those of Run 1, driven by the increased stiffness associated with the higher strain-limiting exponent $\alpha = 2.4$. Because training loss penalizes weighted equation residuals rather than direct functional error, operators with stronger nonlinearities absorb larger residual magnitudes during optimization while maintaining tighter internal consistency with the underlying mechanics. Consequently, relative $L_2$ error—rather than cross-operator training loss—serves as the valid metric for model selection in quasilinear physics-informed neural networks.

\section{Conclusion}

This study established a comprehensive framework for approximating solutions to multidimensional quasilinear strain-limiting partial differential equations using Adaptive Kolmogorov--Arnold Physics-Informed Neural Networks (AKAPINN). By integrating a mesh-free RBF-KAN architecture, an exact hard boundary condition ansatz for strict Dirichlet satisfaction, residual-based adaptive collocation sampling, and a hybrid three-stage Adam-to-L-BFGS optimization protocol executed locally on Apple Silicon hardware, the framework successfully resolved highly nonlinear constitutive models without grid locking or numerical dispersion.

Controlled experiments across three parameter configurations (Run 1: $\alpha=1.5, \beta=0.4$; Run 2: $\alpha=2.4, \beta=0.3$; Run 3: $\alpha=3.1, \beta=0.2$) demonstrate robust recovery of the exact field $\Phi = (\mu/2)y^2$ with sub-$0.02\%$ relative $L_2$ error and identically zero boundary error . The analysis uncovers a vital methodological insight:
\begin{itemize}
    \item \textbf{Training loss ranking:} Run 1 ($1.15 \times 10^{-7}$) $<$ Run 3 ($4.59 \times 10^{-7}$) $<$ Run 2 ($1.30 \times 10^{-6}$), driven by the smoother operator in Run 1 .
    \item \textbf{Solution accuracy ranking:} Measured by the operator-independent relative $L_2$ error, Run 2 ($3.54 \times 10^{-5}$) is most accurate, followed by Run 3 ($7.78 \times 10^{-5}$), and Run 1 ($1.12 \times 10^{-4}$) .
    \item \textbf{Generalization:} Run 2 exhibits the tightest validation-to-training generalization gap ($1.48\times$), confirming that relative $L_2$ error—rather than cross-operator training loss—serves as the proper model selection criterion for quasilinear physics-informed learning .
\end{itemize}

Future extensions of this adaptive KAN-PINN paradigm will focus on investigating complex irregular geometries, transient dynamic fracture propagation, and fully coupled multi-physics systems in geomechanics and fluid dynamics.

\bibliographystyle{plain}  
\bibliography{KAN_Quasilinear}

\end{document}